\PassOptionsToPackage{table}{xcolor}
\documentclass[]{ifacconf}
\usepackage{amssymb,amsmath,graphicx,float}
\usepackage{tikz,pgfplots}
\usetikzlibrary{arrows.meta,shapes.arrows}
\pgfplotsset{compat=1.14}
\usetikzlibrary{arrows,decorations.markings}
\usetikzlibrary{decorations.pathreplacing}
\usetikzlibrary{automata,arrows,positioning,calc}
\usepackage{mathtools}
\usepackage{color} 
\usepackage{tabularx}
\usepackage{algorithm, algpseudocode}

\newtheorem{assumption}{\bf Assumption}

\newtheorem{definition}{\bf Definition}

\newtheorem{problem}{\bf Problem}

\usepackage{times}
\usepackage{lipsum}
\usepackage{nicefrac}
\usepackage{array}
\newcolumntype{P}[1]{>{\centering\arraybackslash}p{#1}}
\allowdisplaybreaks
\usepackage[symbol]{footmisc}
\usepackage{natbib}
\let\ifacconfcaptionwidth\captionwidth
\usepackage[caption=false]{subfig}
\let\captionwidth\ifacconfcaptionwidth

\begin{document}

\begin{frontmatter}
    
    \title{Minimizing the Information Leakage Regarding High-Level Task Specifications}
    
    \author[First]{Michael Hibbard$^1$}
    \author[First]{Yagiz Savas$^1$}
    \author[Second]{Zhe Xu}
    \author[First]{Ufuk Topcu}

    \thanks[]{M. Hibbard and Y. Savas contributed equally to this work.}

    \address[First]{Department of Aerospace Engineering, University of Texas at Austin, TX, USA (email: {\tt\small \{mwhibbard,yagiz.savas,utopcu\}@utexas.edu})}
    \address[Second]{Oden Institute for Computational Engineering and Sciences, University of Texas at Austin, TX, USA. (email: {\tt\small zhexu@utexas.edu})}
    
    \begin{abstract} 
        We consider a scenario in which an autonomous agent carries out a mission in a stochastic environment while passively observed by an adversary. For the agent, minimizing the information leaked to the adversary regarding its high-level specification is critical in creating an informational advantage. 
        We express the specification of the agent as a parametric linear temporal logic formula, measure the information leakage by the adversary's confidence in the agent's mission specification, and propose algorithms to synthesize a policy for the agent which minimizes the information leakage to the adversary. In the scenario considered, the adversary aims to infer the specification of the agent from a set of candidate specifications, each of which has an associated likelihood probability. The agent's objective is to synthesize a policy that maximizes the entropy of the adversary's likelihood distribution while satisfying its specification. We propose two approaches to solve the resulting synthesis problem. The first approach computes the exact satisfaction probabilities for each candidate specification, whereas the second approach utilizes the Fr\'echet inequalities to approximate them. For each approach, we formulate a mixed-integer program with a quasiconcave objective function. We solve the problem using a bisection algorithm. Finally, we compare the performance of both approaches on numerical simulations.
\end{abstract}      


\end{frontmatter}

\section{Introduction}
In environments where privacy and security concerns are of paramount importance, the ability of an agent to deceive an adversary regarding its specification is critical in creating an informational advantage. Pertinent environments for deceptive policies include military operations (\cite{lloyd2003art}), criminal justice (\cite{skolnick1982deception}), and cybersecurity (\cite{carroll2011game}). We explore the concept of deception through the lens of minimizing the information leaked to an adversarial observer regarding the agent's high-level specification.

Specifically, we consider an autonomous agent operating in an environment while passively observed by an adversary. We assume that both the agent and the adversarial observer have knowledge of a set of \textit{specifications}. From this set, the agent maintains a secret \textit{ground-truth specification}, \textit{i.e.}, the specification that the agent actually seeks to satisfy. The adversarial observer attempts to infer the ground-truth specification based on the trajectories of the agent. The agent must behave in such a way as to prevent the adversary from inferring its ground-truth specification. By doing so, the agent may inhibit the adversarial observer from optimally allocating its resources towards preventing the satisfaction of the ground-truth specification.

Consider an autonomous agent that must deliver supplies to one of three possible locations, denoted $A$, $B$, and $C$. The set of candidate specifications for the agent is $\{$\textit{deliver to A}, \textit{deliver to B}, \textit{deliver to C}$\}$. The agent's ground-truth specification is $\{$\textit{deliver to C}$\}$. Although the agent need only travel to $C$ to complete this specification, doing so makes it apparent to an observer that \{\textit{deliver to C}\} is the agent's ground-truth specification. Instead, the agent should travel to each possible location with an equal probability. By doing so, the adversarial observer cannot leverage these probabilities towards inferring which candidate specification is the ground-truth specification.

In this paper, we develop a method for an autonomous agent to synthesize a policy that satisfies the agent's ground-truth specification with desired probability while leading an adversarial observer to infer that each of the candidates is equally likely to be the ground-truth specification. We model the behavior of the agent as a Markov decision process (MDP) (\cite{puterman2014markov}). MDPs are commonly used to model planning and acting in stochastic environments with nondeterministic action selection. Numerous methods exist to synthesize policies for MDPs, which resolve the nondeterminism by prescribing probability distributions for action selection. 

To model the agent's specifications, we use parametric linear temporal logic (pLTL) (\cite{pLTL2014,pLTL}). Standard linear temporal logic (LTL) allows for the formal expression of specifications related to the occurrence of an event, causality between events, and the ordering of successive events (\cite{Baier2008}). pLTL extends LTL by introducing parameterized temporal operators, which allows specifications to be expressed over particular time horizons.

We assume that the adversarial observer assigns a likelihood probability to each candidate according to a simple averaging rule, and set the objective of the agent as maximizing the entropy of the adversary's likelihood distribution. The information-theoretic concept of entropy (\cite{cover2012elements}) measures the average uncertainty of a random variable, an ideal measure for the task at hand. We propose two methods to solve the resulting synthesis problem. The first method exactly computes the probabilities that each specification is completed, which we formulate as a quasiconcave mixed-integer program (MIP) and solve using a bisection method. The exact solution method faces an exponential blow-up in the state space as a function of the number of candidate specifications. The second method we propose avoids the state-space blowup by instead using lower bounds for the probabilities that each specification is satisfied. We formulate this method as an MIP and again solve using a bisection method.

Recently, the works of \cite{savas2018entropy} and \cite{karabag2018least} focused on synthesizing policies that are either unpredictable or difficult for an adversarial observer to infer. These studies focused on the low-level actions rather than on the high-level specifications as we do.
Inferring temporal logic formulas has been extensively studied. For example, \cite{Neider2018} inferred LTL properties classifying a labeled set of trajectories. 
As for the inference of pLTL specifications, \cite{zhe_info} inferred pLTL formulas from a set of trajectories that was \textit{informative} with respect to prior knowledge. Our work is conceptually similar to these studies; however, we seek to make the inference problem as difficult as possible. 


\section{Preliminaries}\label{sec_prelim}
\noindent \textbf{Notation.} We denote the set $\{1,2,\ldots\}$ of natural numbers and the set $(-\infty,\infty)$ of real numbers by $\mathbb{N}$ and $\mathbb{R}$, respectively. For a given logical formula, $\top$ and $\bot$ denote that the formula is true and false, respectively. For a set $S$, we denote its power set by $2^{S}$. Finally, for $N$$\in$$\mathbb{N}$, we denote the set $\{1,2,\ldots,N\}$ by $[N]$. 

\subsection{Markov Decision Processes}

{\setlength{\parindent}{0cm}
\begin{definition}
A \textit{Markov decision process} (MDP) is defined by the tuple $\mathcal{M}$$=$$( S, s_0, \mathcal{A}, \mathcal{P}, \mathcal{AP}, \mathcal{L})$ where $S$ is a finite set of states, $\mathcal{A}$ is a finite set of actions, $s_0$ is a unique initial state, $\mathcal{P}$$:$$ S$$\times$$\mathcal{A}$$\times$$S$$\rightarrow$$[0,1]$ is a transition function such that $\sum_{s'\in S}\mathcal{P}(s,a,s')$$=$$1$ for all $a$$\in$$\mathcal{A}$ and $s$$\in$$S$, $\mathcal{AP}$ is a set of atomic propositions, and $\mathcal{L}$$:$$S $$\rightarrow$$2^{\mathcal{AP}}$ is a labeling function.
\end{definition}}
We denote the transition probability $\mathcal{P}(s,a,s')$ by $\mathcal{P}_{s,a,s'}$. The size of an MDP is the number of triples $(s,a,s')$$\in$$S\times \mathcal{A}\times S$ in which $\mathcal{P}_{s,a,s'}$$>$$0$.

{\setlength{\parindent}{0cm}
\begin{definition}
A \textit{policy} $\pi$ for an MDP $\mathcal{M}$ is a sequence $\pi$$=$$(d_1,d_2,d_3,\ldots)$ where each $d_t$$:$$S$$\times$$\mathcal{A}$$\rightarrow$$[0,1]$ is a mapping such that $\sum_{a\in\mathcal{A}}d_t(s,a)$$=$$1$ for all $s$$\in$$S$. For an MDP $\mathcal{M}$, we denote the set of all admissible policies by $\Pi(\mathcal{M})$.
\end{definition}}
A \textit{stationary} policy satisfies $\pi$$=$$(d_1,d_1,d_1,\ldots)$. We denote the probability of choosing an action $a$$\in$$\mathcal{A}$ in a state $s$$\in$$S$ under a stationary policy $\pi$ by $\pi(s,a)$. 

For an arbitrary length $L$$\in$$\mathbb{N}$, we refer to a sequence of states $\varrho^{\pi}$$:=$$s_0s_1\ldots s_L$ generated in $\mathcal{M}$ under a policy $\pi$$\in$$\Pi(\mathcal{M})$ as a \textit{trajectory}, which starts from the initial state $s_0$ and satisfies $\sum_{a_t\in\mathcal{A}}d_t(s_t,a_t)\mathcal{P}_{s_t,a_t,s_{t+1}}$$>$$0$ for all $0$$\leq$$t$$<$$L$.

\subsection{Parametric Linear Temporal Logic}
Following \cite{pLTL2014}, the syntax of parametric linear temporal logic (pLTL) is defined recursively as
\[
\begin{split}
\phi:=&\top\mid p \mid\lnot\phi\mid\phi_{1}\wedge\phi_{2}\mid \bigcirc\phi\mid \phi_1\mathcal{U}\phi_2 \mid \lozenge_{\sim i} \phi,
\label{syntax}
\end{split}
\]
where $p$ is an atomic proposition, $\lnot$ and $\wedge$ stand for negation and conjunction, respectively, $\bigcirc$ and $\mathcal{U}$ are temporal operators representing \textquotedblleft next\textquotedblright~and \textquotedblleft until\textquotedblright, respectively, $\Diamond_{\sim i}$ is a parameterized temporal operator representing `` \textit{parameterized eventually}", where $\sim$$\in$$\{\ge,\le\}$, and $i$$\in$$\mathbb{N}$ is a temporal parameter. We recursively define the logical connective $\vee$ (disjunction), and temporal operators $\Diamond$ (eventually), $\Box$ (always), $\Box_{\sim i}$ (parameterized always) and $\mathcal{U}_{\sim i}$ (parameterized until) from the aforementioned operators (\cite{pLTL2014}). Furthermore, for $i_1$$<$$i_2$, we define the parameterized temporal operators $\Diamond_{[i_1, i_2]}$ and $\Box_{[i_1, i_2]}$ such that, for a formula $\phi$,
$\Diamond_{[i_1, i_2]}\phi=\Diamond_{\ge i_1}\phi\wedge\Diamond_{\le i_2}\phi$ and 
$\Box_{[i_1, i_2]}\phi=\Box_{\ge i_1}\phi\wedge\Box_{\le i_2}\phi$.

For an MDP $\mathcal{M}$ under a policy $\pi$$\in$$\Pi(\mathcal{M})$, a trajectory $\varrho^{\pi}$$=$$s_0s_1\ldots s_L$ generates a word $w^{\pi}$$:=$$w_0w_1\ldots w_L$, where $w_k$$=$$\mathcal{L}(s_k)$ for all $0\leq$$k$$\leq$$L$. For a pLTL formula $\phi$ and a trajectory $\varrho^{\pi}$ at time index $k$$\leq$$L$, the satisfaction relation $(w^{\pi},k)\models\phi$ is defined recursively as  
\[
\begin{split}
(w^{\pi},k)\models p \quad\mbox{iff}\quad& p \in\mathcal{L}(s_k),\\
(w^{\pi},k)\models\lnot\phi\quad\mbox{iff}\quad & (w^{\pi},k)\not\models\phi,\\
(w^{\pi},k)\models\phi_{1}\wedge\phi_{2}\quad\mbox{iff}\quad & (w^{\pi},k)\models\phi_{1}\quad\\& \mbox{and}\quad(w^{\pi},k)\models\phi_{2},\\
(w^{\pi},k)\models\bigcirc\phi\quad\mbox{iff}\quad& (w^{\pi},k+1)\models\phi,\\
(w^{\pi},k)\models\phi_{1}\mathcal{U}\phi_{2}\quad\mbox{iff}\quad &  \exists
 k'\ge k,\ (w^{\pi},k')\models\phi_{2},\\
  & \text{and}\  \forall k^{\prime\prime}\in[k, k'], (w^{\pi},k^{\prime\prime})\models\phi_{1}, \\
(w^{\pi},k)\models\Diamond_{\sim i}\phi\quad\mbox{iff}\quad & \exists
k'\sim k+i, ~(w^{\pi},k')\models\phi.
\end{split}
\] 
If the satisfaction relations are evaluated at time index $k=0$, then we simply write $w^{\pi}\models\phi$. For an scpLTL formula $\phi$, the set $\{\varrho^{\pi}$$:$$w^{\pi}$$\models$$\phi\}$ is measurable (\cite{Baier2008}). We denote $\text{Pr}^{\pi}_{\mathcal{M}}(w$$\models$$\phi)$ as the probability that a word $w$, generated by an MDP $\mathcal{M}$ under a policy $\pi$$\in$$\Pi(\mathcal{M})$, satisfies a pLTL formula $\phi$; \textit{i.e.}, $w$$\in$$\{\varrho^{\pi}: w^{\pi} \models \phi\}$.

As discussed in \cite{zhe_info}, \textit{syntactically co-safe pLTL} (scpLTL) formulas are a special class of pLTL formulas that can be satisfied by words of finite length. The syntax of scpLTL is defined recursively as
	\[
	\begin{split}
	\phi:=&\top\mid\pi\mid\lnot\pi\mid\phi_{1}\wedge\phi_{2}\mid\phi_{1}\vee
	\phi_{2}\mid \bigcirc\phi\mid \Diamond\phi\mid \phi_1\mathcal{U}\phi_2\\& \mid \Diamond_{\sim i}\phi\mid \Box_{\le i}\phi\mid \phi_1\mathcal{U}_{\sim i}\phi_2.
	\end{split}
	\]
Because scpLTL is a restriction of pLTL, the satisfaction relation of scpLTL formulas can be derived from the satisfaction relation of general pLTL formulas.

\section{Problem Formulation}
We consider an agent operating in a stochastic environment whose behavior is modeled by an MDP. The agent aims to complete a task, expressed as a \textit{ground-truth} scpLTL specification $\phi^{\star}$, with desired probability $\Gamma$$\in$$(0,1)$, while in the presence of an adversarial observer. The adversary aims to infer the task of the agent through observations of its trajectory. Aware of the adversary's objective, the agent aims to complete its task with the desired probability while simultaneously minimizing the information leaked to the adversary about the task.

Let $\phi:=\{\phi_{1},\phi_{2},\ldots,\phi_{N} \}$ be a set of scpLTL specifications such that $\phi^{\star}$$\in$$\phi$. The adversary has a finite set $\phi_{can}$$\subseteq$$\phi$ of \textit{candidate scpLTL specifications}, which it uses to describe the task of the agent. In particular, let $\beta$$\in$$(0,1)$ be a constant candidacy threshold, and $\pi$$\in$$\Pi(\mathcal{M})$ be the agent's policy. An scpLTL specification $\phi_i$$\in$$\phi$ is a candidate, i.e., $\phi_i$$\in$$\phi_{can}$, if and only if $\text{Pr}^{\pi}_{\mathcal{M}}(w $$\models$$\phi_i)$$\geq$$\beta$. In other words, a specification is a candidate if the trajectories followed by the agent satisfy the specification with at least probability $\beta$.

We \textit{assume} that, to each candidate specification $\phi_i$$\in$$\phi_{can}$, the adversary assigns a likelihood probability
\begin{align}\label{prob_assign}
  &\text{Pr}(\phi_i=\phi^{\star} | \phi_i\in \phi_{can}):=\nonumber \\
  &\qquad \qquad \qquad \quad \frac{\text{Pr}^{\pi}_{\mathcal{M}}(s_0\models \phi_i)\mathbb{I}\{\phi_i \in \phi_{can}\}}{\sum_{\phi_i\in \phi}\text{Pr}^{\pi}_{\mathcal{M}}(s_0\models \phi_i)\mathbb{I}\{\phi_i \in \phi_{can}\}}
\end{align}
where $\mathbb{I}\{s$$\in$$S \}$ is an indicator function such that $\mathbb{I}\{s$$\in$$S \}$$:=$$1$ if $s$$\in$$S$ and $\mathbb{I}\{s$$\in$$S \}$$:=$$0$ otherwise. The probability assignment \eqref{prob_assign} is a simple averaging rule representing the adversary's confidence in the candidate being the ground-truth specification. The adversary may also measure its confidence level using a distribution different from \eqref{prob_assign}; \textit{e.g.}, a Boltzmann distribution. In that case, the solution techniques introduced in this paper can still be utilized to synthesize a policy minimizing the adversary's information about the task. However, depending on the distribution, the synthesis of such a policy may require one to employ different computational methods.

We use the adversary's certainty on ground-truth specification as the measure of the information leakage. For a given policy $\pi$$\in$$\Pi(\mathcal{M})$, let $\text{Pr}_{\mathcal{M},\pi, \phi_i}$$:=$$\text{Pr}(\phi_i$$=$$\phi^{\star} | \phi_i$$\in$$\phi_{can})$. We measure the adversary's \textit{uncertainty} on the specification $\phi^{\star}$ by the entropy
\begin{align}
    H^{\pi}(\phi_{can}):=-\sum_{\phi_i\in\phi_{can}} \text{Pr}_{\mathcal{M},\pi, \phi_i}\log \text{Pr}_{\mathcal{M},\pi, \phi_i}
\end{align}
of the distribution $\text{Pr}(\phi_i$$=$$\phi^{\star} | \phi_i$$\in$$\phi_{can})$. The rationale behind this choice can be better understood by recalling that the entropy of a random event is the lower bound on the average number of bits required to describe the outcomes of the event (\cite{cover2012elements}). Moreover, this lower bound is maximized when the probability distribution associated with the event is uniform. By following a policy maximizing $H^{\pi}(\phi_{can})$, the agent satisfies all candidate specifications $\phi_i$$\in$$\phi_{can}$ with \textit{similar} probabilities, making it more difficult for the adversary to guess the ground-truth specification $\phi^{\star}$ with high confidence. 

\begin{problem} \label{entropyProb}
 	Given an MDP $\mathcal{M}$, a set of candidate scpLTL formulas $\phi$$=$$\{\phi_{1},\phi_{2},\ldots,\phi_{N} \}$, a ground-truth formula $\phi^{\ast}$$\in$$\phi$, and constants $\Gamma,\beta$$\in$$(0,1)$ such that $\Gamma$$\geq$$\beta$, synthesize a policy $\pi$$\in$$\Pi(\mathcal{M})$ that solves the following problem:
	\begin{subequations}
	\begin{align}\label{obj_prob_1}
	    \underset{\pi\in\Pi(\mathcal{M})}{\text{maximize}} \ \ &H^{\pi}(\phi_{can})\\
	    \text{subject to:} \ \ & \text{Pr}^{\pi}_{\mathcal{M}}(w\models \phi^{\star})\geq \Gamma\\
	    & \phi_i \in \phi_{can} \ \ \iff \ \ \text{Pr}^{\pi}_{\mathcal{M}}(w\models \phi_i)\geq \beta \label{cons_prob_1}
	\end{align}
	\end{subequations}
\end{problem} 
~ 
Intuitively, a policy solving the problem defined in \eqref{obj_prob_1}-\eqref{cons_prob_1} maximizes the uncertainty of the adversary about the agent's task while ensuring that the agent completes the task with desired probability.

\section{An Exact Solution Method}\label{section_exact}

We now present an exact solution method for the problem defined in \eqref{obj_prob_1}-\eqref{cons_prob_1}. First, we construct a product MDP on which the satisfaction probability of each specification $\phi_i$$\in$$\phi$ can be verified. We then formulate a nonlinear optimization problem on this product MDP, whose solution provides a policy solving the problem defined in \eqref{obj_prob_1}-\eqref{cons_prob_1}.

\subsection{Product MDP} \label{section_product_mdp}
We construct a product MDP
in three steps. First, we construct a deterministic finite automaton for each specification $\phi_i$. Second, we form an expanded MDP whose state labels track the stage number of the underlying process. Finally, we take the product of the expanded MDP with each of the automata constructed in the first step.

For any scpLTL specification $\phi_i$ with fixed parameters, one can construct a deterministic finite automaton with the input alphabet $2^{\mathcal{AP}}$ which accepts a word $w^{\pi}$ if and only if (iff) $w^{\pi}$ satisfies the specification $\phi_i$, $i.e.$, $w^{\pi}$$\models$$\phi_i$ (\cite{KupfermanVardi2001}).

{\setlength{\parindent}{0cm}
\begin{definition}
A \textit{deterministic finite automaton} (DFA) is a tuple $A$$=$$(Q, q_{0},$$2^{\mathcal{AP}},$$\delta,$$\mathcal{F})$, where $Q$ is a finite set of states, $q_{0}$ is a unique initial state, $2^{\mathcal{AP}}$ is an alphabet, $\delta:$$Q$$\times$$2^{\mathcal{AP}}$$\rightarrow$$Q$ is a transition function, and $\mathcal{F}$$\subseteq$${Q}$ is a finite set of accepting states. 
\label{DFA}
\end{definition} }

For a given scpLTL formula $\phi_i$$\in$$\phi$, we denote its corresponding DFA by $A_{\phi_i}$. Without loss of generality (w.l.o.g.), we assume that the accepting states $\mathcal{F}$ of $A_{\phi_i}$ are absorbing, \textit{i.e.}, $\delta(q,p)$$=$$q$ for all $q$$\in$$\mathcal{F}$ and $p$$\in$$2^{\mathcal{AP}}$.
We do not lose generality since an input word is accepted by a DFA $A_{\phi_i}$ iff it has a finite prefix that reaches an accepting state on $A_{\phi_i}$. The continuation of the word after that prefix has no effect on its acceptance by $A_{\phi_i}$.

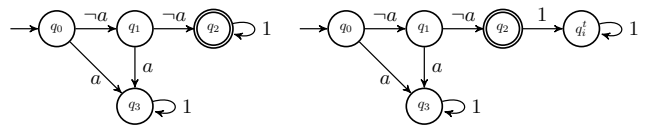
\begin{figure}[b!]\centering
\scalebox{0.77}{\begin{tikzpicture}[->, >=stealth', auto, semithick, node distance=2cm]

    \tikzstyle{every state}=[fill=white,draw=black,thick,text=black,scale=0.7]

    \node[state,initial,initial text=] (s_0) {$q_0$};
    \node[state] (s_1) [right=7mm of s_0]  {$q_1$};
     \node[state,accepting] (s_2) [right=7mm of s_1]  {$q_2$};
     \node[state] (s_3) [below=7mm of s_1]  {$q_3$};

\path
(s_0)  edge  node[]{$\lnot a$}     (s_1)
(s_0)  edge  node[below]{$a$}     (s_3)
(s_1)  edge      node[]{$a$}     (s_3)
(s_1)  edge   node{$\lnot a$}     (s_2)
(s_2)  edge  [loop right=10]    node{$1$}     (s_2)
(s_3)  edge  [loop right=10]    node{$1$}     (s_3);
\end{tikzpicture}}
\scalebox{0.77}{\begin{tikzpicture}[->, >=stealth', auto, semithick, node distance=2cm]

    \tikzstyle{every state}=[fill=white,draw=black,thick,text=black,scale=0.7]

    \node[state,initial,initial text=] (s_0) {$q_0$};
    \node[state] (s_1) [right=7mm of s_0]  {$q_1$};
     \node[state,accepting] (s_2) [right=7 mm of s_1]  {$q_2$};
     \node[state] (s_3) [below=7mm of s_1]  {$q_3$};
     \node[state] (s_4) [right=7 mm of s_2]  {$q_i^t$};

\path
(s_0)  edge  node[]{$\lnot a$}     (s_1)
(s_0)  edge  node[below]{$a$}     (s_3)
(s_1)  edge      node[]{$a$}     (s_3)
(s_1)  edge   node{$\lnot a$}     (s_2)
(s_2)  edge    node{$1$}     (s_4)
(s_3)  edge  [loop right=10]    node{$1$}     (s_3)
(s_4)  edge  [loop right=10]    node{$1$}     (s_4);
\end{tikzpicture}}
\caption{ An example construction of the modified DFA $\overline{A}_{\phi_i}$ for the scpLTL specification $\phi_i$$=$$\square_{\leq 2}\lnot a$. (Left) The nominal DFA $A_{\phi}$ where $q_2$$\in$$\mathcal{F}$ is the only accepting state. (Right) The modified DFA $\overline{A}_{\phi}$. }\label{modified_auto_ex}
\end{figure}

We modify a given DFA $A_{\phi_i}$ by augmenting $Q$ with a terminal state $q_i^t$ which is absorbing and reachable only from the states in $\mathcal{F}$. Specifically, the modified DFA $\overline{A}_{\phi_i}$ has the finite set of states $\overline{Q}$$:=$$Q\cup\{q_i^t\}$, with a transition function $\overline{\delta}$$:$$Q\times 2^{\mathcal{AP}}$$\rightarrow$$Q$ defined by
\begin{align*}
    \overline{\delta}(q,p):=\begin{cases} q_i^t & \text{if}\ q\in \mathcal{F}\cup \{q_i^t\}\\
    \delta(q,p) & \text{otherwise}.
    \end{cases}
\end{align*}
In Fig. \ref{modified_auto_ex}, we provide an example construction of the modified DFA $\overline{A}_{\phi_i}$ for the scpLTL formula $\phi_i$$=$$\square_{\leq 2} \lnot a$ where $\{a\}$$\in$$\Sigma$. 

We now form the expanded MDP whose state labels tracks the stage number of the underlying process so that the satisfaction of a given scpLTL specification can be verified.

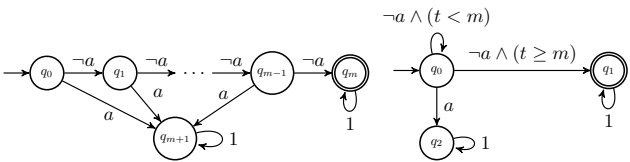
\begin{figure}[b!]\centering
\scalebox{0.72}{\begin{tikzpicture}[->, >=stealth', auto, semithick, node distance=2cm]

    \tikzstyle{every state}=[fill=white,draw=black,thick,text=black,scale=0.7]

    \node[state,initial,initial text=] (s_0) {$q_0$};
    \node[state] (s_1) [right=7mm of s_0]  {$q_1$};
    \node[] (s_11) [right=7mm of s_1]  {$\ldots$};
     \node[state] (s_2) [right=7mm of s_11]  {$q_{m-1}$};
      \node[state,accepting] (s_3) [right=7mm of s_2]  {$q_{m}$};
     \node[state] (s_4) [below right =7mm and 5 mm of s_1]  {$q_{m+1}$};

\path
(s_0)  edge  node[]{$\lnot a$}     (s_1)
(s_1)  edge  node[]{$\lnot a$}     (s_11)
(s_11)  edge  node[]{$\lnot a$}     (s_2)
(s_0)  edge  node[below]{$a$}     (s_4)
(s_1)  edge      node[]{$a$}     (s_4)
(s_2)  edge      node[above]{$a$}     (s_4)
(s_2)  edge  node[]{$\lnot a$}    (s_3)
(s_3)  edge  [loop below=10]    node{$1$}     (s_3)
(s_4)  edge  [loop right=10]    node{$1$}     (s_4);
\end{tikzpicture}}
\scalebox{0.72}{\begin{tikzpicture}[->, >=stealth', auto, semithick, node distance=2cm]

    \tikzstyle{every state}=[fill=white,draw=black,thick,text=black,scale=0.7]

    \node[state,initial,initial text=] (s_0) {$q_0$};
    \node[state,accepting] (s_1) [right=25mm of s_0]  {$q_1$};
     \node[state] (s_2) [below=7 mm of s_0]  {$q_2$};

\path
(s_0)  edge  [loop above=10]    node{$\lnot a \land (t<m)$}     (s_0)
(s_0)  edge  node[]{$\lnot a \land (t \geq m)$}     (s_1)
(s_0)  edge  node[]{$ a$}     (s_2)
(s_2)  edge  [loop right=10]    node{$1$}     (s_2)
(s_1)  edge  [loop below=10]    node{$1$}     (s_1);
\end{tikzpicture}}
\caption{The effect of the input alphabet on the size of the automaton for the scpLTL specification $\phi_i$$=$$\square_{\leq m} \lnot a$ where $m$$\in$$\mathbb{N}$ is a constant. (Left) The DFA has the input alphabet $2^{\{a\}}$. (Right) The DFA has the input alphabet $2^{\{a\}\cup [m]}$. }\label{reduction_auto_ex}
\end{figure}

{\setlength{\parindent}{0cm}
\begin{definition} \label{expanded_MDP_def} Let $\mathcal{M}$$=$$(S,s_0,\mathcal{A},\mathcal{P},\mathcal{AP},\mathcal{L})$ be an MDP and $[\mathcal{T}]$$:=$$\{1,2,\ldots,\mathcal{T}\}$ be an index set. The \textit{expanded MDP} $\mathcal{M}$$\times$$[T]$$=$$(S^{[\mathcal{T}]}, s_0^{[\mathcal{T}]}, \mathcal{A}, \mathcal{P}^{[\mathcal{T}]}, \mathcal{L}^{[\mathcal{T}]},\mathcal{AP}^{[\mathcal{T}]})$ is a tuple where $S^{[\mathcal{T}]}=S\times [\mathcal{T}]$, $s_0^{[\mathcal{T}]}=(s_0,1)$, $\mathcal{P}^{[\mathcal{T}]}((s,t),a,(s',t'))=$
\begin{equation*}
    \begin{cases}\mathcal{P}_{s,a,s'} & \text{if}\  t<\mathcal{T} \ \land\  t'=t+1\\
    \mathcal{P}_{s,a,s'} & \text{if}\  t=\mathcal{T} \ \land \ t'=t \\
    0 & \text{otherwise,}\end{cases}
\end{equation*}
$\mathcal{L}^{[\mathcal{T}]}((s,t))=\mathcal{L}(s)\cup\{t\}$, and $\mathcal{AP}^{[\mathcal{T}]}=\mathcal{AP}\cup [\mathcal{T}]$.
\end{definition}}
We note that by choosing $2^{\mathcal{AP}^{[\mathcal{\mathcal{T}}_i]}}$ instead of $2^{\mathcal{AP}}$ as the input alphabet $\Sigma$, one can dramatically decrease the number of states in the DFA $A_{\phi_i}$. In Fig. \ref{reduction_auto_ex}, we demonstrate the significance of the input alphabet on the size of a DFA corresponding to the formula $\phi_i$$=$$\square_{\leq m} \lnot a$. To reduce the size of the state-space in the solution of the problem \eqref{obj_prob_1}-\eqref{cons_prob_1}, we verify the satisfaction of a given formula $\phi_i$ over $\mathcal{M}$$\times$$[\mathcal{\mathcal{T}}_i]$ instead of $\mathcal{M}$.

To verify if the probability that the trajectories followed by an agent on $\mathcal{M}$$\times$$[\mathcal{T}_i]$ satisfies a specification $\phi_i$ exceeds a desired threshold, one can construct a product MDP and verify whether the agent's trajectories reach the accepting states of the product MDP with desired probability. Note in the following definition that we abuse the notation for the expanded MDP $\mathcal{M}$$\times$$[\mathcal{T}_i]$.

{\setlength{\parindent}{0cm}
\begin{definition}\label{product_def}
Let $\mathcal{M}$$\times$$[\mathcal{T}_i]$$=$$(S,s_0,\mathcal{A},\mathcal{P},\mathcal{AP},\mathcal{L})$ be an expanded MDP and $\overline{A}_{\phi_i}$$=$$(\overline{Q}_i,q^i_0,2^{\mathcal{AP}},\overline{\delta}_i, \mathcal{F}_i)$ be a modified DFA. The \textit{product MDP} $\mathcal{M}$$\times$$[\mathcal{T}_i]$$\times$$\overline{A}_{\phi_i}$$=$$(S_p, s_{0_p}, \mathcal{A}, \mathbb{P}, \mathcal{AP}$, $\mathcal{L}_p,$$\mathcal{F}_p)$ is a tuple where $S_p$$=$$S $$\times$$ \overline{Q}_i$, $s_{0_p}=(s_0,q)$ such that $q=\overline{\delta}_i(q^i_0,\mathcal{L}(s_0))$, $\mathbb{P}((s,q), a, (s',q'))$$=$
\begin{equation*}
    \begin{cases} \mathcal{P}_{s,a,s'} & \text{if} \quad q'=\overline{\delta}_i(q,\mathcal{L}(s')) \\ 0 & \text{otherwise}, \end{cases},
\end{equation*}
$\mathcal{L}_p((s,q))=\{q\}$, and $\mathcal{F}_p$$=$$S\times \mathcal{F}_i$.
\end{definition}}

A product MDP $\mathcal{M}\times[\mathcal{T}_i]\times \overline{A}_{\phi_i}$ may contain states that are not reachable from the initial state. Unreachable states have no effect in the analysis of the MDP. These states can be found in time polynomial in the size of $\mathcal{M}\times [\mathcal{T}_i]\times \overline{A}_{\phi_i}$ by graph search algorithms, \textit{e.g.}, breadth-first search, and can subsequently be removed from the product MDP w.l.o.g. We hereafter assume that there is no unreachable state in $\mathcal{M}\times [\mathcal{T}_i] \times \overline{A}_{\phi_i}$.

For a given specification $\phi_i$ with a fixed parameter set ${\bf{p}}_i$, let $\mathcal{T}_i$$:=$$\max {\bf{p}}_i$ be the maximum element of ${\bf{p}}_i$, \textit{e.g.}, ${\bf{p}}_i$$=$$\{4,8\}$ and $\mathcal{T}_i$$=$$8$ for $\phi_i$$=$$\square_{\leq 4} a \land \lozenge_{\geq 8} b$. We note that for nested formulas, the parameter set can be defined recursively. As an example, for $\phi_i$$=$$\lozenge_{[a,b] }\square_{[c,d]} p$, letting $\phi_j$$:=$$\square_{[c,d]} p$, we have ${\bf{p}}_j$$=$$\{c,d\}$, and  ${\bf{p}}_i$$=$$\{a+c, a+d, b+c, b+d\}$. For an MDP $\mathcal{M}$ and a set $\phi$$=$$\{\phi_1,\phi_2,\ldots,\phi_N\}$ of specifications, we form the product MDP $\mathcal{M}_p$$:=$$\mathcal{M}$$\times$$[\mathcal{T}]$$\times$$\overline{A}_{\phi_1}$$\times$$\overline{A}_{\phi_2}$$\times$$\ldots$$\times$$\overline{A}_{\phi_N}$ by recursively applying Definition \ref{product_def}, where $\mathcal{T}$$:=$$\max_{i\in[N]}\mathcal{T}_i$. In this construction, the input alphabet to each DFA $\overline{A}_{\phi_i}$ is $2^{\mathcal{AP}^{[\mathcal{T}]}}$.
\subsection{Policy Synthesis: An Optimization Problem}\label{opt_problem_section}
After constructing the product MDP $\mathcal{M}_p$ on which the satisfaction probability of each specification $\phi_i$ can be verified, we now provide a nonlinear optimization problem whose solution provides a policy solving the problem defined in \eqref{obj_prob_1}-\eqref{cons_prob_1}. 

Let the tuple $\bf{s}$$:=$$(s,t,q_1,q_2,\ldots, q_N)$ denote a state in $\mathcal{M}_p$ such that $s$$\in$$S$, $t$$\in$$[\mathcal{T}]$, and $q_i$$\in$$Q_i$ for all $i$$\in$$\{1,2,\ldots,N\}$. We denote the $k^{th}$ element of the tuple ${\bf{s}}$ by ${\bf{s}}[k]$, \textit{e.g.}, ${\bf{s}}[1]$$=$$s$, ${\bf{s}}[2]$$=$$t$, ${\bf{s}}[3]$$=$$q_1$ and ${\bf{s}}[N+2]$$=$$q_N$. Moreover, with an abuse of notation, we denote the transition function of $\mathcal{M}_p$ by $\mathbb{P}$. 

We partition the states of $\mathcal{M}_p$ into the disjoint sets $B$ and $S_p\backslash B$. Let $B$ be the set of states ${\bf{s}}$$\in$$S_p$ such that 
\begin{align}
    \sum_{\substack{\bf{s}' \in S_p: \\
    {\bf{s}}'[i+2]={\bf{s}}[i+2]\  \forall i\in[N]}} \mathbb{P}_{{\bf{s}},a,{\bf{s}}'}=1 \label{abs_check}
\end{align}
for all $a$$\in$$\mathcal{A}$. The set $B$ is a collection of states $\bf{s}$ whose elements ${\bf{s}}[k$$+$$2]$ correspond to the automata states $q_k$ that are absorbing. Once a state ${\bf{s}}$$\in$$B$ is reached by the agent, we know that each specification $\phi_i$ is either satisfied or violated by the agent. Note that the set $B$ can be computed in time polynomial in the size of $\mathcal{M}_p$, as condition \eqref{abs_check} can be verified by simply checking whether the automata elements of a state are absorbing or not.  

We assume w.l.o.g. that $\phi_1$$=$$\phi^{\star}$, \textit{i.e.}, $\phi_1$ is the ground-truth specification $\phi^{\star}$. Let $\alpha$$:$$S_p$$\rightarrow$$[0,1]$ be a function such that $\alpha(s_{0_p})$$=$$1$ and $\alpha(s_{0_p})$$=$$0$ otherwise, \textit{i.e.}, $\alpha$ is the initial state distribution of $\mathcal{M}_p$. The optimization problem is given as:
\begin{subequations}
\begin{align}\label{opt_obj}
    & \underset{\lambda({\bf{s}},a), x(i)}{\text{maximize}} -\sum_{i\in[N]} \frac{\nu(i)}{\sum_{j\in [N]}\nu(j)}\log \Bigg(\frac{\nu(i)}{\sum_{j\in[N]}\nu(j)}\Bigg) \raisetag{24pt}\\
    & \text{subject to:} \nonumber \\ \label{opt_cons_flow}
    & \forall {\bf{s}}\in S_p\backslash B, \  \sum_{a\in\mathcal{A}}\lambda({\bf{s}},a) -\sum_{{\bf{s}}'\in S_p}\sum_{a\in\mathcal{A}} \mathbb{P}_{{\bf{s}}',a,{\bf{s}}} \lambda({\bf{s}}',a)=\alpha({\bf{s}}) \  \raisetag{10pt}\\ \label{opt_cons_mu}
    &\forall i\in [N],\ \mu(i)=\sum_{{\bf{s}}\in S_p:\  {\bf{s}}[i+2]\in  \mathcal{F}_i}\ \sum_{a\in\mathcal{A}}\lambda({\bf{s}},a) \ \qquad \qquad  \raisetag{24pt}\\ \label{opt_cons_ground_truth}
    & \qquad \qquad \qquad \quad \ \mu(3)\geq \Gamma \\ \label{opt_cons_beta}
    &\forall i\in [N],\ \qquad \quad \mu(i)\geq \beta x(i)  \\ \label{opt_cons_nu}
   &\forall i\in [N], \ \qquad \quad \nu(i)=\mu(i)x(i)\\ \label{opt_cons_lambda_ineq}
   & \forall {\bf{s}}\in S_p, \forall a\in\mathcal{A}, \lambda({\bf{s}},a)\geq 0\\  \label{opt_cons_x_binary}
   &\forall i\in [N],\  \qquad \quad x(i)\in \{0,1\}
\end{align}
\end{subequations}

The decision variables in the above optimization problem are $\lambda(\textbf{s},a)$ for each ${\bf{s}}$$\in$$S_p$ and $a$$\in$$\mathcal{A}$, and $x(i)$ for each $i$$\in$$[N]$. The variables $\mu(i)$ and $\nu(i)$ are functions of $\lambda({\bf{s}},a)$ and $x(i)$, defined in \eqref{opt_cons_mu} and \eqref{opt_cons_nu}, respectively, to simplify the notation. The variable  $\lambda({\bf{s}},a)$ corresponds to the expected number of times the state-action pair $({\bf{s}},a)$ is visited (\cite{Marta}). In particular, we have the relation
\begin{align}
   \lambda({\bf{s}},a)=\sum_{t=1}^{\infty}\text{Pr}^{\pi}_{\mathcal{M}_p}(S_t={\bf{s}}, A_t=a| S_1=s_{0_p})
\end{align}
where the policy $\pi$$\in$$\Pi(\mathcal{M}_p)$ is defined as
\begin{align}\label{opt_policy_construct}
    \pi({\bf{s}},a):=\begin{cases}\frac{\lambda({\bf{s}},a)}{\sum_{a'\in \mathcal{A}}\lambda({\bf{s}},a')} & \text{if} \ \sum_{a'\in \mathcal{A}}\lambda({\bf{s}},a')>0\\
    1/ \lvert \mathcal{A}\rvert & \text{otherwise}.
    \end{cases}
\end{align}
For more details on the variable $\lambda({\bf{s}},a)$, we refer the reader to \cite[Chapter 6]{puterman2014markov}, \cite[Chapter 2]{altman}, and \cite{Marta}. Finally, the binary variable $x(i)$, under the constraints \eqref{opt_cons_flow}-\eqref{opt_cons_x_binary}, satisfies the relation 
\begin{align}
    x(i)=\begin{cases} 1 & \text{if} \ \text{Pr}^{\pi}_{\mathcal{M}}(w\models\phi_i)\geq\beta\\
    0 & \text{otherwise}.
    \end{cases}
\end{align}

Constraint \eqref{opt_cons_flow} is traditionally referred to as the ``flow constraint" (\cite{Marta}), which ensures that the number of times the agent leaves a state is equal to the number of times it enters that state. Constraint \eqref{opt_cons_mu} defines the variable $\mu(i)$, the probability of reaching an accepting state of the automaton $\overline{A}_{\phi_i}$. Although the variable $\lambda({\bf{s}},a)$ is the expected number of visits to $({\bf{s}},a)$, because we use the modified automaton $\overline{A}_{\phi_i}$ in the product MDP $\mathcal{M}_p$, $\lambda({\bf{s}},a)$ corresponds to the reachability probability for states ${\bf{s}}$ satisfying ${\bf{s}}[i+2]$$\in$$\mathcal{F}_i$. Constraints \eqref{opt_cons_ground_truth} and  \eqref{opt_cons_beta} ensure, respectively, that the ground-truth specification $\phi_1$ is satisfied by at least probability $\Gamma$, and that if $x(i)$$=$$1$, we have $\phi_i$$\in$$\phi_{can}$. Constraint \eqref{opt_cons_nu} defines the variable $\nu(i)$, which is equal to the satisfaction probability of the specification $\phi_i$ if $x(i)$$=$$1$ and zero otherwise. Finally, constraints \eqref{opt_cons_lambda_ineq} and \eqref{opt_cons_x_binary} define the feasible domains of the decision variables.

The objective function \eqref{opt_obj} is the entropy of the probability distribution $\nu(i)/\sum_{j\in [N]}\nu(j)$, which, under the constraints \eqref{opt_cons_flow}-\eqref{opt_cons_x_binary}, is equal to the right hand side of \eqref{prob_assign}. Specifically, it can be seen from the constraints \eqref{opt_cons_mu}-\eqref{opt_cons_nu} that we have $\nu(i)$$=$$\text{Pr}^{\pi}_{\mathcal{M}}(w\models \phi_i)\mathbb{I}\{\phi_i$$\in $$\phi_{can}\}$.

We note that, once an optimal solution $\lambda^{\star}({\bf{s}},a)$ to the problem \eqref{opt_obj}-\eqref{opt_cons_x_binary} is computed, one can obtain an optimal policy $\pi^{\star}$$\in$$\Pi(\mathcal{M}_p)$ on the product MDP $\mathcal{M}_p$ using the construction given in \eqref{opt_policy_construct}. Then, using the one-to-one correspondence between the policies on $\mathcal{M}$ and $\mathcal{M}_p$ (see, \textit{e.g.}, \cite{Baier2008},\cite{Wolff2012}), we can finally construct a policy on $\mathcal{M}$, which solves the problem \eqref{obj_prob_1}-\eqref{cons_prob_1}.

\subsection{Policy Synthesis: A Solution Approach}\label{section_opt_solution}
The nonlinear optimization problem \eqref{opt_obj}-\eqref{opt_cons_x_binary} has a special structure which can be exploited to utilize off-the-shelf optimization toolboxes such as GUROBI (\cite{gurobi}) and MOSEK (\cite{mosek}) for obtaining a global optimal solution. In this section, we provide an algorithm, based on a bisection method \cite[Chapter 4]{boyd2004convex}, that allows the utilization of such toolboxes.

We begin with the exact relaxation of the constraint \eqref{opt_cons_nu}. Note that  \eqref{opt_cons_nu} is a bilinear constraint since both $\mu(i)$ and $x(i)$ are variables in the optimization problem. Such constraints are not handled by most off-the-shelf toolboxes. However, recalling that $\mu(i)$ represents the probability of reaching an accepting state of the automaton $\overline{A}_{\phi_i}$, we know that $0\leq \mu(i)\leq 1$. Using this additional information, we can replace each constraint \eqref{opt_cons_nu}, with its corresponding McCormick envelope (\cite{McCormick1976}), given by the following inequalities
\begin{align}\label{mccormick1}
   & \nu(i)\geq 0,&&\quad \nu(i)\leq x(i), \\ \label{mccormick2}
    & \nu(i)\leq \mu(i),&& \quad \nu(i)\geq x(i)+\mu(i)-1.
\end{align}
Note that, using the above inequalities, we have $\nu(i)$$=$$0$ if $x(i)$$=$$0$, and $\nu(i)$$=$$\mu(i)$ if $x(i)$$=$$1$. Therefore, the relaxation of the constraint \eqref{opt_cons_nu} with the above inequalities is exact. Moreover, since the above constraints are affine in the variables $\mu(i)$ and $x(i)$, they can now be handled by off-the-shelf toolboxes.

Next, we utilize the quasiconcavity of the objective function in \eqref{opt_obj} in the variables $\nu(i)$.
\begin{figure}[b!]\centering\scalebox{0.32}{
\includegraphics[]{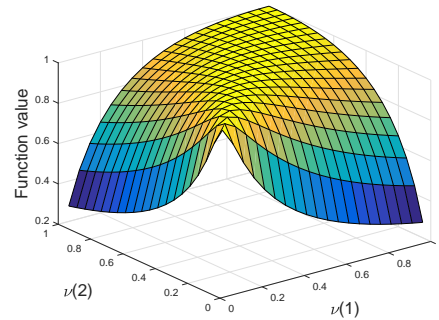}}
\caption{The function $f(\boldsymbol{\nu})$ is not concave. However, $f(\boldsymbol{\nu})$ has a different useful property, which is quasiconcavity.  }
\label{quasi_concave}
\end{figure}
A quasiconcave function is formally defined below. For additional details on convex sets and functions, we refer the reader to \cite{boyd2004convex}.
{\setlength{\parindent}{0cm}
\begin{definition} (\cite{boyd2004convex})\label{quasi_con_def} A function $g$$:$$\mathcal{D}$$\rightarrow$$\mathbb{R}$ is called \textit{quasiconcave} if its domain $\mathcal{D}$ and all its superlevel sets $\Omega_{\theta}$$:=$$\{x\in \mathcal{D}: g(x)\geq \theta\}$ for $\theta$$\in$$\mathbb{R}$ are convex.
\end{definition}}

Let $\boldsymbol{\nu}$$:=$$[\nu(1),\nu(2),\ldots,\nu(N)]$ be a vector of variables $\nu(i)$ for $i\in[N]$, and $f$$:$$\mathbb{R}^N_+\rightarrow \mathbb{R}$ be a function such that
\begin{align*}
    f(\boldsymbol{\nu}):=-\sum_{i\in[N]} \frac{\nu(i)}{\sum_{j\in [N]}\nu(j)}\log \Bigg(\frac{\nu(i)}{\sum_{j\in[N]}\nu(j)}\Bigg).
\end{align*}
We emphasize that the function $f(\boldsymbol{\nu})$ is not concave, as illustrated in Fig. \ref{quasi_concave} for $\boldsymbol{\nu}$$\in$$[0,1]^2$. By defining functions $f_1$$:$$\mathbb{R}^N_+$$\rightarrow$$\mathbb{R}$ and $f_2$$:$$\mathbb{R}^N_+$$\rightarrow$$\mathbb{R}$ such that
\begin{align*}
    f_1(\boldsymbol{\nu})&:=-\sum_{i\in[N]}\nu(i)\log \Bigg(\frac{\nu(i)}{\sum_{j\in[N]}\nu(j)}\Bigg),\\
    f_2(\boldsymbol{\nu})&:=\sum_{j\in [N]}\nu(j),
\end{align*}
we obtain the relation $ f(\boldsymbol{\nu})$$=$$f_1(\boldsymbol{\nu})/ f_2(\boldsymbol{\nu})$. Convexity of the sublevel set $\Omega_{\theta}$$:=$$\{\boldsymbol{\nu}$$\in$$\mathbb{R}^N_+$$:$$ f(\boldsymbol{\nu})$$\geq$$\theta\}$ for any $\theta$$\in$$\mathbb{R}$ follows from the fact that $f_1(\boldsymbol{\nu})$$\geq$$\theta$$ f_2(\boldsymbol{\nu})$ defines a convex region since the functions $f_1(\boldsymbol{\nu})$ and $f_2(\boldsymbol{\nu})$ are, respectively, concave and affine over their domains (\cite{boyd2004convex}). We thus conclude the quasiconcavity of $f(\boldsymbol{\nu})$ from Definition \ref{quasi_con_def}. 

We are now ready to introduce an iterative algorithm for the solution of \eqref{opt_obj}-\eqref{opt_cons_x_binary}, which is a variant of the bisection method for quasiconcave optimization (Algorithm 4.1 in \cite{boyd2004convex}).
Let $v^{\star}$ be the optimal value of the problem in \eqref{opt_obj}-\eqref{opt_cons_x_binary}, and $u$$\in$$\mathbb{R}$ be an arbitrarily large constant which satisfies $u$$\geq$$v^{\star}$. Moreover, let $l$$\in$$\mathbb{R}$ be a constant such that $l$$\leq$$v^{\star}$, \textit{e.g.}, $l$$=$$0$. At each iteration of the algorithm, we set $\theta$$:=$$(u+l)/2$ and solve the feasibility problem given in \eqref{feasibility_prob}. If the problem has a feasible solution, in the next iteration of the algorithm, we set $u$$:=$$\theta$, otherwise, we set $l$$:=$$\theta$. The algorithm terminates when the stop condition $u-l$$\leq$$\epsilon$ is satisfied, where $\epsilon$$>$$0$ is a constant tolerance parameter.
\begin{equation}\label{feasibility_prob}
\begin{aligned}
     \text{find} \qquad  &\theta \\
     \text{subject to:} \quad &f_1(\boldsymbol{\nu})\geq \theta f_2(\boldsymbol{\nu})\\
     &\boldsymbol{\nu}=[\nu(1),\nu(2),\ldots,\nu(N)]\\
    & \eqref{opt_cons_flow}, \eqref{opt_cons_mu}, \eqref{opt_cons_ground_truth}, \eqref{opt_cons_beta}, \eqref{opt_cons_lambda_ineq}, \eqref{opt_cons_x_binary}, \eqref{mccormick1},\eqref{mccormick2}
\end{aligned}
\end{equation}

As mentioned in Section \ref{opt_problem_section}, an optimal policy $\pi^{\star}$$\in$$\Pi(\mathcal{M}_p)$ on the product MDP $\mathcal{M}_p$ can be obtained using the construction given in \eqref{opt_policy_construct}, once an optimal solution $\lambda^{\star}({\bf{s}},a)$ to the problem \eqref{opt_obj}-\eqref{opt_cons_x_binary} is computed using the algorithm defined in (\ref{feasibility_prob}).

\section{An Approximate Solution Method}
Although the method presented in Section \ref{section_exact} provides an exact solution to the Problem 1, it requires one to form the product MDP, which is the product of the expanded MDP with the automata corresponding to each specification $\phi_i$. The construction of the product MDP is, in general, a computationally expensive operation; therefore, for practical purposes, it is desirable to develop algorithms that synthesize policies directly on the expanded MDP.
In this section, focusing on a subset of scpLTL specifications, we present a method that conservatively approximates the satisfaction probabilities of the specifications $\phi_i$ and allows one to synthesize policies on the expanded MDP.  

Throughout this section, we restrict our attention to a subset of scpLTL specifications with the following assumption.

{\setlength{\parindent}{0cm}
\begin{assumption}\label{assumption_1}
Each scpLTL specification $\phi_{i}$$\in$$\phi$ we consider has one of four possible forms: $\lozenge_{[a,b]}p$, $\square_{[a,b]}p$, $\lozenge_{[a,b]} \square_{[c,d]}p$, or $\square_{[a,b]} \lozenge_{[c,d]}p$, where $a,b,c,d$$\in$$\mathbb{N}$ and $p$$\in$$\mathcal{AP}$.
\end{assumption}}

We now present the Fr\'echet inequalities (\cite{frechet,hailperin1965best}), which allow us to conservatively approximate the satisfaction probability of a specification $\phi_i$. For each $i$$\in$$[K]$ where $K$$\in$$\mathbb{N}$, let $E_i$ be a logical proposition, and $e_i$ be the probability that the proposition $E_i$ is true. Then,
\begin{align}\label{frechet_1}
    \text{Pr}\Big(\bigwedge_{i=1}^K E_i \Big)&\geq \max\Big\{0, \sum_{i=1}^K e_i-(n-1)\Big\},\\ \label{frechet_2}
    \text{Pr}\Big(\bigvee_{i=1}^K E_i \Big)&\geq \max\{e_1,e_2,\ldots,e_K\}.
\end{align}

These lower bounds are the best possible bounds if nothing is known about the events $E_1,E_2,\ldots,E_K$ except that their probabilities are $e_1,e_2,\ldots,e_K$, respectively (\cite{hailperin1965best}). A remarkable property of these lower bounds is that they are in terms of the satisfaction probabilities $e_i$ of the subformulas $E_i$ only. If one can represent a logical formula $E$ as a conjunction or disjunction of the subformulas $E_i$ whose satisfaction probability can be computed easily, then by ensuring that the lower bound exceeds a desired threshold $\beta$, one can guarantee that the formula $E$ is satisfied with at least probability $\beta$. 

In what follows, we form an expanded MDP on which, instead of measuring the satisfaction probability of a specification $\phi_i$, we measure the satisfaction probabilities of subformulas of $\phi_i$ corresponding to \textit{each time step}. As an example, instead of measuring the satisfaction probability of $\phi_i$$=$$\square_{[1,3]}p$, we measure the probability that the predicate $p$ holds true at a given time step $1$$\leq$$t$$\leq$$3$. Then, using the syntax of scpLTL specifications, we utilize these measurements to derive the lower bound on the satisfaction probability of $\phi_i$.

Recall from Section \ref{section_product_mdp} that for a given specification $\phi_i$ with a fixed parameter set ${\bf{p}}_i$, we have $\mathcal{T}_i$$=$$\max {\bf{p}}_i$. For an MDP $\mathcal{M}$ and a set $\phi$$=$$\{\phi_1,\phi_2,\ldots,\phi_N\}$ of specifications, using Definition \ref{expanded_MDP_def}, we form the expanded MDP $\overline{\mathcal{M}}$$:=$$\mathcal{M}$$\times$$[\mathcal{T}$$+$$1]$ where $\mathcal{T}$$:=$$\max_{i\in[N]}\mathcal{T}_i$. On $\overline{\mathcal{M}}$, we can measure the satisfaction probability of a predicate $p$$\in$$\mathcal{AP}$ at time $t$$\in$$[\mathcal{T}]$ by the expected number of visits to states ${\bf{s}}$$\in$$S^{[\mathcal{T}+1]}$ such that ${\bf{s}}[2]$$=$$t$ and $p$$\in$$\mathcal{L}({\bf{s}})$. 

To synthesize a policy $\pi$$\in$$\Pi(\overline{\mathcal{M}})$ on the expanded MDP $\overline{\mathcal{M}}$, we solve a modified version of the problem \eqref{opt_obj}-\eqref{opt_cons_x_binary} on the state-space $S^{[\mathcal{T}+1]}$ of $\overline{\mathcal{M}}$. In particular, for each specification $\phi_i$, we replace the corresponding constraint \eqref{opt_cons_mu} with a series of other constraints. Recall that the variable $\mu(i)$ in \eqref{opt_cons_mu} is equal to the probability of satisfying the specification $\phi_i$. Instead of using the exact satisfaction probability, for each specification form in Assumption \ref{assumption_1}, we introduce a set of constraints which ensure that $\mu(i)$ is a lower bound on the actual satisfaction probability. 

$\bullet \, $$\phi_{i}$$=$$\lozenge_{[k_{1},k_{2}]}p$. We first introduce the variables $\eta(t)$$\in$$\mathbb{R}$ for each $t$$\in$$\mathbb{N}$ such that $k_1$$\leq$$t$$\leq$$k_2$. Using the syntax of pLTL specifications, we can show that the lower bound in \eqref{frechet_2} is equal to $\mu(i)$, defined by the following constraints:
\begin{subequations}
\begin{align}\label{reach_exact_first}
     &\sum_{\substack{\bf{s} \in S^{[\mathcal{T}+1]}: \\ {\bf{s}}[2]=t, \ p\in \mathcal{L}({\bf{s}})}} \, \, \sum_{a \in A} \lambda({\bf{s}},a) = \eta(t), \\ \label{mu_max_equality}
     &\mu(i)=\max\Big\{\eta(k_1),\eta(k_1+1),\ldots,\eta(k_2)\Big\}.
\end{align}
\end{subequations}
\noindent
In the above constraints, each variable $\eta(t)$ captures the probability that the formula $\phi_{i}$ holds at the particular time step $t$. To utilize the off-the-shelf toolboxes for encoding the above constraints, we need to relax the constraint \eqref{mu_max_equality}. We do so by replacing \eqref{mu_max_equality} with the following set of constraints
\begin{subequations}
\begin{align}\label{reach_relax_first}
    &\forall k_1\leq t\leq k_2,\  \mu(i)\geq \eta(t);  \quad \mu(i)=\sum_{t=k_1}^{k_2}y(t)\eta(t), \\ 
    & \forall k_1\leq t\leq k_2,\   y(t)\in\{0,1\}; \quad \sum_{t=k_1}^{k_2}y(t)=1.
\end{align}
\end{subequations}
The above relaxation is exact. In \eqref{reach_relax_first}, the term $y(t)\eta(t)$ is bilinear as both $y(t)$ and $\eta(t)$ are variables. However, since we know that $0$$\leq$$\eta(t)$$\leq$$1$ from \eqref{reach_exact_first}, by defining an extra variable $\gamma(t)$$:=$$y(t)\eta(t)$, we can represent each term $\gamma(t)$ exactly with its corresponding McCormick envelope given in \eqref{mccormick1}-\eqref{mccormick2}. 

$\bullet \, $$\phi_{i}$$=$$\square_{[k_{1},k_{2}]}p$. We first introduce the variables $\eta(t)$$\in$$\mathbb{R}$ for each $t$$\in$$\mathbb{N}$ such that $k_1$$\leq$$t$$\leq$$k_2$. Using the syntax of pLTL specifications, we can show that the lower bound in \eqref{frechet_1} is equal to $\mu(i)$, defined by the following constraints:
\begin{subequations}
\begin{align}
     &\sum_{\substack{\bf{s} \in S^{[\mathcal{T}+1]}: \\ {\bf{s}}[2]=t, \ p\in \mathcal{L}({\bf{s}})}} \, \, \sum_{a \in A} \lambda({\bf{s}},a) = \eta(t), \\  \label{mu_min_ineq}
     &\mu(i)=\max\Bigg\{0, \sum_{t=k_1}^{k_2}\eta(t)-(k_2-k_1)\Bigg\}.
\end{align}
\end{subequations}
To utilize the off-the-shelf toolboxes for encoding the above constraints, we need to relax the constraint \eqref{mu_min_ineq}. We do so by replacing \eqref{mu_min_ineq} with the following set of constraints

\begin{subequations}
\begin{align}
    & \mu(i)\geq 0; \ \ \quad \mu(i)\geq \Bigg(\sum_{t=k_1}^{k_2}\eta(t)-(k_2-k_1)\Bigg), \\ 
    & y\in\{0,1\}; \quad \mu(i)= y\Bigg(\sum_{t=k_1}^{k_2}\eta(t)-(k_2-k_1)\Bigg).
\end{align}
\end{subequations}
The above relaxation is exact, but it involves the bilinear terms $y\eta(t)$. However, once again, by introducing new variables $\gamma(t)$$:=$$y\eta(t)$, we can represent each term $\gamma(t)$ exactly with its corresponding McCormick envelope given in \eqref{mccormick1}-\eqref{mccormick2}. 

$\bullet \, $$\phi_{i}$$=$$\lozenge_{[k_{1},k_{2}]}$$\square_{[k_{3},k_{4}]}p$. We first introduce the variables $\eta(t)$$\in$$\mathbb{R}$ for each $t$$\in$$\mathbb{N}$ such that $k_1$$+$$k_3$$\leq$$t$$\leq$$k_2$$+$$k_4$, and $\zeta(m)$ for each $m$$\in$$\mathbb{N}$ such that $k_1$$\leq$$m$$\leq$$k_2$. Using the syntax of pLTL specifications and both of the bounds in \eqref{frechet_1}-\eqref{frechet_2}, we can obtain a lower bound $\mu(i)$ on the satisfaction probability of $\phi_i$ using the following constraints:

\begin{subequations}
\begin{align}
     &\sum_{\substack{\bf{s} \in S^{[\mathcal{T}+1]}: \\ {\bf{s}}[2]=t, \ p\in \mathcal{L}({\bf{s}})}} \, \, \sum_{a \in A} \lambda({\bf{s}},a) = \eta(t), \\  \label{similar1}
     &\zeta(m)=\max\Bigg\{0, \sum_{t=m+k_3}^{m+k_4}\eta(t)-(k_4-k_3)\Bigg\},\\ \label{similar2}
     &\mu(i)=\max\Big\{\zeta(k_1),\zeta(k_1+1),\ldots,\zeta(k_2)\Big\}.
\end{align}
\end{subequations}
We can perform the relaxation of the constraints in \eqref{similar1}-\eqref{similar2} by introducing new binary variables and subsequently using the corresponding McCormick envelopes as previously explained in the relaxation of the specifications $\lozenge_{[k_1,k_2]}p$ and $\square_{[k_1,k_2]}p$.

$\bullet \ \phi_{i}$$=$$\square_{[k_{1},k_{2}]}$$\lozenge_{[k_{3},k_{4}]}p$. We first introduce the variables $\eta(t)$$\in$$\mathbb{R}$ for each $t$$\in$$\mathbb{N}$ such that $k_1+k_3$$\leq$$t$$\leq$$k_2+k_4$, and $\zeta(m)$ for each $m$$\in$$\mathbb{N}$ such that $k_1$$\leq$$m$$\leq$$k_2$. Using the syntax of pLTL specifications and both of the bounds in \eqref{frechet_1}-\eqref{frechet_2}, we can obtain a lower bound $\mu(i)$ on the satisfaction probability of $\phi_i$ using the following constraints:
\begin{subequations}
\begin{align}
     &\sum_{\substack{\bf{s} \in S^{[\mathcal{T}+1]}: \\ {\bf{s}}[2]=t, \ p\in \mathcal{L}({\bf{s}})}} \, \, \sum_{a \in A} \lambda({\bf{s}},a) = \eta(t), \\  \label{similar3}
     &\zeta(m)=\max\Big\{\eta(m+k_3),\eta(m+k_3+1),\ldots, \eta(m+k_4)\Big\},\\ \label{similar4}
     &\mu(i)=\max\Bigg\{0, \sum_{m=k_1}^{k_2}\zeta(m)-(k_2-k_1)\Bigg\}.
\end{align}
\end{subequations}
Again, we perform the relaxation of the constraints in \eqref{similar3}-\eqref{similar4} by introducing new binary variables and using the corresponding McCormick envelopes as explained in the relaxation of the specifications $\lozenge_{[k_1,k_2]}p$ and $\square_{[k_1,k_2]}p$.

Finally, after replacing each constraint \eqref{opt_cons_mu} in \eqref{opt_obj} with its corresponding set of constraints introduced in this section, we solve the resulting nonlinear optimization problem using the bisection method presented in Section \ref{section_opt_solution}.

\section{Numerical Examples}

We now provide several examples to demonstrate the efficacy of the proposed solution methods. For each example, we use a tolerance of $\epsilon$$=$$1$$\times$$10^{-4}$ for the bisection method. We use the GUROBI solver with the CVX (\cite{cvx}) interface to solve the exact and approximate optimization problems.

\begin{figure}[b!]
     \centering
     \scalebox{0.8}{
     \subfloat{\begin{tikzpicture}[every node/.style={minimum size=.6cm-5*\pgflinewidth, outer sep=0pt}]
            
            \draw[step=0.6cm,color=black] (0,0) grid (3.6,3.6);
        
            \node[fill=white] at (0.3,3.3) {$\tiny{I}$};
            \node[fill=blue] at (0.3,0.3) {};
            \node[fill=blue] at (0.9,0.3) {};
            \node[fill=blue] at (0.3,0.9) {};
            \node[fill=blue] at (0.9,0.9) {};
            \node[fill=red] at (2.7,2.7) {};
            \node[fill=red] at (2.7,3.3) {};
            \node[fill=red] at (3.3,2.7) {};
            \node[fill=red] at (3.3,3.3) {};
            \node[fill=yellow] at (1.5,1.5) {};
            \node[fill=yellow] at (1.5,2.1) {};
            \node[fill=yellow] at (2.1,1.5) {};
            \node[fill=yellow] at (2.1,2.1) {};
            \node[fill=green] at (2.7,0.3) {};
            \node[fill=green] at (2.7,0.9) {};
            \node[fill=green] at (3.3,0.3) {};
            \node[fill=green] at (3.3,0.9) {};

     \end{tikzpicture}}
     \label{example1grid}}
     \hspace{6mm}
     \scalebox{0.8}{
     \subfloat{\begin{tikzpicture}[->, >=stealth', auto, semithick, node distance=2cm]

        \tikzstyle{every state}=[fill=white,draw=black,thick,text=black,scale=0.7]

        \node[state] (s1) {$S_{1}$};
        \node[state] (s2) [above left = 7mm and 7mm of s1]  {$S_{2}$};
        \node[state] (s3) [above right = 7 mm and 7mm of s1]  {$S_{3}$};
        \node[state] (s4) [below right = 7mm and 7mm of s1]  {$S_{4}$};
        \node[state] (s5) [below left = 7mm and 7mm of s1]  {$S_{5}$};
        \node[] (dummy) [below = 12.25mm of s1] {};
        
        \path
        (s1)    edge    (s2)
        (s2)    edge    (s1)
        (s1)    edge    (s3)
        (s3)    edge    (s1)
        (s1)    edge    (s4)
        (s4)    edge    (s1)
        (s1)    edge    (s4)
        (s1)    edge    (s5)
        (s5)    edge    (s1)
        (s2)    edge    (s3)
        (s3)    edge    (s2)
        (s3)    edge    (s4)
        (s4)    edge    (s3)
        (s4)    edge    (s5)
        (s5)    edge    (s4)
        (s5)    edge    (s2)
        (s2)    edge    (s5)
        (s1)    edge    [loop below=10] (s1)
        (s2)    edge    [loop left=10] (s2)
        (s3)    edge    [loop right=10] (s3)
        (s4)    edge    [loop right=10] (s4)
        (s5)    edge    [loop left=10] (s5);
        
        \end{tikzpicture}}
     \label{example2network}}
     \caption{Environments considered in examples. (Left) Gridworld considered in the resupply mission. (Right) MDP considered in the surveillance mission.}
     \label{steady_state}
\end{figure}
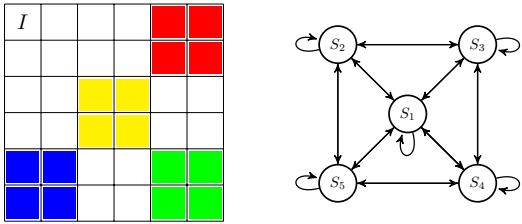

\begin{table}[b!]
    \begin{center}
    \caption{Sets of specifications used for the examples. ${\phi}^{*}$ indicates the ground-truth specification.} \label{table:specComp}
    \rowcolors{2}{gray!25}{white}
    \begin{tabular}{|c|p{0.2\textwidth}|} 
            \hline
            Example & Specifications \\
            \hline
            Resupply-1 & $\phi^{*}$$:$$\square_{[9,10]}blue$ \newline $\phi_{2}$$:$$\square_{[29,30]}$$red$ \\
            \hline
            Resupply-2 & $\phi^{*}$$:$$\square_{[9,10]}blue$ \newline $\phi_{2}$$:$$\square_{[16,18]}yellow$ \newline
            $\phi_{3}$$:$$\square_{[23,25]}green$ \newline
            $\phi_{4}$$:$$\square_{[29,30]}$$red$ \\
            \hline
            Surveillance & $\phi^{*}$$:$$\square_{[1,10]}$$\lozenge_{[0,5]}$$blue$ \newline
            $\phi_{2}$$:$$\square_{[1,10]}$$\lozenge_{[0,5]}$$red$\newline
            $\phi_{3}$$:$$\square_{[1,10]}$$\lozenge_{[0,5]}$$yellow$ \newline
            $\phi_{4}$$:$$\square_{[1,10]}$$\lozenge_{[0,5]}$$green$\\
            \hline
    \end{tabular}
    \end{center}
\end{table}

\subsection{A Resupply Mission}

We first consider an autonomous agent operating on the gridworld shown in Figure \ref{example1grid}. The colored states represent different bases that the agent can travel to. The agent's mission is to resupply the blue base, which we encode as the pLTL formula ``$\square_{[a,b]}blue$"; i.e., the agent must reach the blue base at a specified time $a$ and remain there until its supplies are unloaded after $b$$-$$a$ time steps. Due to the presence of an adversarial observer, the agent must additionally attempt to obfuscate which base it actually delivers the supplies to. By doing so, the adversary is least able to infer which base actually received the supplies.

The agent is assumed to start in the upper left corner of the gridworld in state $I$. In each state, the agent can select one of four possible actions: move left, move right, move up, or move down. Once the agent has selected an action, it transitions to its desired state with probability 0.99, while slipping to the left, to the right, or backwards each with probability $\nicefrac{.01}{3}$.

We study two cases for the specifications of the agent. We first consider that the agent only seeks to prevent information leakage about its ground-truth specification by additionally visiting the red base. We subsequently consider that the agent seeks to prevent information leakage by additionally visiting the green and yellow bases as well. For each set of specifications, we let $\Gamma$$=$$0.95$ and $\beta$$=$$0.8$. The sets of specifications are provided in Table \ref{table:specComp}, under ``Resupply-1" and ``Resupply-2", respectively. For each set of specifications, we run the exact and approximate methods to solve the optimization problem (\ref{opt_obj})-(\ref{opt_cons_x_binary}).

\begin{table*}[]\label{TableResults}
    \centering
    \caption{Number of continuous and binary variables, solution times, probabilities of satisfaction, resulting entropy, and the size of the set $\phi_{can}$ for exact and approximate solution methods.}\label{table:resultsComp}
    \rowcolors{2}{gray!25}{white}
    \begin{tabular}{ |c| P{0.07\textwidth}| P{0.07\textwidth}| P{0.037\textwidth}| P{0.033\textwidth}| P{0.075\textwidth}| P{0.085\textwidth}| P{0.08\textwidth}| P{0.059\textwidth}| P{0.059\textwidth}| P{0.034\textwidth}| P{0.034\textwidth}|}
            \hline
            \scriptsize{Example} & \scriptsize{Num. of Var.} \newline \scriptsize{exact} & \scriptsize{Num. of Var.} \newline \scriptsize{approx.} & \scriptsize{Time} \newline \scriptsize{exact} & \scriptsize{Time} \newline \scriptsize{approx.} & \scriptsize{$\text{Pr}^{\pi}_{\mathcal{M}}$$(w$$\models$$\phi^{*})$} \newline \scriptsize{exact} & \scriptsize{$\text{Pr}^{\pi}_{\mathcal{M}}$$(w$$\models$$\phi^{*})$} \newline \scriptsize{approx., comp.} & \scriptsize{$\text{Pr}^{\pi}_{\mathcal{M}}$$(w$$\models$$\phi^{*})$} \newline \scriptsize{approx., actual} & \scriptsize{$H^{\pi}(\phi_{can})$} \newline \scriptsize{exact} & \scriptsize{$H^{\pi}(\phi_{can})$} \newline \scriptsize{approx.} & \scriptsize{$|\phi_{can}|$} \newline \scriptsize{exact} & \scriptsize{$|\phi_{can}|$} \newline \scriptsize{approx.} \\
            \hline
             \scriptsize{Resupply-1} & \scriptsize{4870 con.} \newline \scriptsize{1 binary} & \scriptsize{3079 con.} \newline \scriptsize{1 binary} & \scriptsize{9.25s} & \scriptsize{6.75s} & \scriptsize{0.950} & \scriptsize{0.950} & \scriptsize{0.971} & \scriptsize{1.000} & \scriptsize{0.999} & \scriptsize{2} & \scriptsize{2} \\
            \hline
             \scriptsize{Resupply-2} & \scriptsize{6980 con.} \newline \scriptsize{3 binary} & \scriptsize{3125 con.} \newline \scriptsize{3 binary} & \scriptsize{53.71s} & \scriptsize{25.11s} & \scriptsize{0.950} & \scriptsize{0.950} & \scriptsize{0.971} & 1.999 & \scriptsize{1.999} & \scriptsize{4} & \scriptsize{4} \\
            \hline
            \scriptsize{Surveillance} & \scriptsize{22209 con.} \newline \scriptsize{3 binary} & \scriptsize{617 con.} \newline \scriptsize{239 binary} & \scriptsize{148.37s} & \scriptsize{26.59s} & \scriptsize{0.950} & \scriptsize{0.951} & \scriptsize{0.991} & \scriptsize{1.999} & \scriptsize{1.999} & \scriptsize{4} & \scriptsize{4} \\
            \hline
    \end{tabular}
\end{table*}

Table \ref{table:resultsComp} lists the relevant output information for each solution method and specification set. We note that the number of variables in the optimization problem is after GUROBI completed presolving the problem. Because of how the approximation for specifications of the form $\square_{[a,b]}$ was constructed, the two solution methods for each specification set have the same number of binary variables. However, as the approximate solution method does not require taking the product with each specification automaton, the number of continuous variables in its corresponding optimization problem scales better than that of the exact solution method and requires less time to solve. 

The approximate solution method performs nearly as well as the exact solution method at minimizing the information leakage about the ground-truth specification. Both solution methods obtain the maximum size of the candidate set $\phi_{can}$ and nearly obtain the maximum-entropy upper bounds of 1 and 2 bits for each of the two solution methods, respectively. 

\subsection{A Surveillance Mission}

We now consider an agent that must repeatedly surveil an outpost containing sensitive information on the boundary of its base. Specifically, the agent operates on the MDP shown in Figure \ref{example2network}, where four outposts surround the central base. We assume that the adversarial observer does not know which of the outposts contains the sensitive information. For this reason, the agent must additionally surveil the three non-sensitive outposts. By doing so, the adversary cannot use the fact that the agent visits an outpost towards inferring which outpost contains the sensitive information. Thus, the adversary cannot optimally allocate its resources towards infiltrating the correct outpost.

We assume that the agent's initial state is in the central state $S_{1}$. In each state, the agent can either remain in its current state or transition to a neighboring state, where it transitions with probability 1. We again set $\Gamma$$=$$0.95$ and $\beta$$=$$0.8$, respectively. We use the pLTL structure ``$\square_{[a,b]}\lozenge_{[c,d]}outpost_{i}$" to encode the surveillance specification; i.e., at each time step over the time horizon, the agent must eventually visit an outpost within $[c,d]$ time steps. The specifications for the agent are listed in Table \ref{table:specComp} under ``Surveillance". For this set of specifications, we again run the exact and approximate methods to solve the optimization problem (\ref{opt_obj})-(\ref{opt_cons_x_binary}).

Table \ref{table:resultsComp} shows the comparison of the output between the exact and approximate solution methods for the set of surveillance specifications. Although the approximate solution method uses a large number of binary variables compared to the exact solution method, it is still able to solve the optimization problem (\ref{opt_obj})-(\ref{opt_cons_x_binary}) much quicker than the exact solution method is able to. The approximate and exact solution methods perform similarly well in minimizing the information leakage about which of the outposts contained sensitive information. Both solution methods obtain the maximum number of elements in the candidate set $\phi_{can}$ and nearly achieve the upper bound on the maximum entropy of 2 bits.

\section{Conclusions}
We study the problem of synthesizing a policy for an autonomous agent that leaks the minimum amount of information regarding its high-level task specification to an adversarial observer. We measure the information leakage as the adversary's confidence that a candidate mission specification is the ground-truth mission specification. Modelling the inference problem of the adversary as an averaging rule, we formulate the problem of the agent as a mixed-integer program with a quasiconcave objective function, and develop two methods for its solution. The first method exactly computes the probabilities that a specification is satisfied by the agent, whereas the second method approximates these probabilities using the Fr\'echet inequalities. We provide two numerical examples to demonstrate the efficacy of the proposed solution methods in minimizing the information leakage.

\bibliography{floodLTL}     

\end{document}